\documentclass[12pt,reqno]{amsart}

\usepackage{times}
\usepackage[T1]{fontenc}
\usepackage{mathrsfs}
\usepackage{latexsym}
\usepackage[dvips]{graphics}
\usepackage[dvips]{graphicx}
\usepackage{epsfig}
\usepackage{amsmath,amsfonts,amsthm,amssymb,amscd}
\input amssym.def
\input amssym.tex
\usepackage{color}

\newcommand\be{\begin{equation}}
\newcommand\ee{\end{equation}}
\newcommand\bea{\begin{eqnarray}}
\newcommand\eea{\end{eqnarray}}
\newcommand\nbea{\begin{eqnarray*}}
\newcommand\neea{\end{eqnarray*}}
\newcommand\bi{\begin{itemize}}
\newcommand\ei{\end{itemize}}
\newcommand\ben{\begin{enumerate}}
\newcommand\een{\end{enumerate}}

\newcommand{\Z}{\ensuremath{\mathbb{Z}}}

\numberwithin{equation}{section}

\newcommand{\bal}{\begin{align}}
\newcommand{\eal}{\end{align}}

\begin{document}

\title{If a prime divides a product....}

\author{Steven J. Miller}\email{Steven.J.Miller@williams.edu}
\address{Department of Mathematics and Statistics, Williams College,
Williamstown, MA 01267}

\author{Cesar E. Silva}\email{csilva@williams.edu}
\address{Department of Mathematics and Statistics, Williams College,
Williamstown, MA 01267}

\subjclass[2010]{ (primary),  (secondary)}

\keywords{}

\date{\today}

\thanks{We thank Jared Hallet, for asking a good question in Silva's analysis class, which led to this work. The first named author was partly supported by NSF grant DMS0970067. }

\begin{abstract} One of the greatest difficulties encountered by all in their first proof intensive class is subtly assuming an unproven fact in a proof. The purpose of this note is to describe a specific instance where this can occur, namely in results related to unique factorization and the concept of the greatest common divisor.
\end{abstract}

\maketitle

%%%%%%%%%%%%%%%%%%%%%%%%%%%%%%%%%%%%%%%%%%%%%%%%%%%%%%%%%%%%%%%%%%%%%%%%%%%%%%%%%%%%%%%%%%%%%%%%%%%%%%%%%%%%

%%%%%%%%%%%%%%%%%%%%%%%%%%%%%%%%%%%%%%%%%%%%%%%%%%%%%%%%%%%%%%%%%%%%%%%%%%%%%%%%%%%%%%%%%%%%%%%%%%%%%%%%%%%%

The Fundamental Theorem of Arithmetic states that every integer exceeding 1 can be written uniquely as a product of prime powers. There are may ways to prove this important theorem, and many applications. One of the most important consequences of unique factorization is in studying the Riemann zeta function, defined by $$\zeta(s) \ = \ \sum_{n=1}^\infty \frac1{n^s}$$ if the real part of $s$ exceeds 1. If this is your first time seeing this function, it shouldn't be clear at all why it is worth studying, or why the Clay Mathematics Institute \cite{Cl} is offering one million dollars for a proof about the location of its zeros!\footnote{Although $\zeta(s)$ is initially only defined when the real part of $s$ is greater than 1, we can use complex analysis to analytically continue the function to be defined for all complex $s \neq 1$, and this extension agrees with the original series expansion when the real part of $s$ exceeds 1. This extended function is known to vanish at the negative even integers. The celebrated Riemann Hypothesis asserts that the only other zeros of the continued function have real part equal to 1/2.} This turns out to be one of the most important functions in number theory. The reason is that unique factorization implies that we may rewrite the zeta function as a product: if the real part of $s$ exceeds 1 then we also have $$\zeta(s) \ = \ \prod_{p\ {\rm prime}} \left(1 - \frac1{p^s}\right)^{-1};$$ see Chapter 3 of \cite{MT-B}. This relationship is the starting point of many investigations. The reason is that the integers are built up from the primes, and the distribution and properties of the primes are difficult and mysterious. The integers, however, are quite well understood. For example, it is not immediately apparent what  the next prime after the prime 17483 is (it is 17489); however, it is quite easy to find the next integer after 17483! The equivalence of the sum and the product allow us to transfer information from the integers to information of the primes.

The point of all of this is that unique factorization is a very important property, one which needs to be carefully proved.
This note grew out of some conversations between the two authors over passages in the second named author's upcoming book on real analysis \cite{Si}. Specifically, we saw how easy it was to implicitly assume results in arguments and proofs. The purpose here is to discuss not just the proof of unique factorization and similar results, but also to highlight how easy it is to subtly assume a fact when trying to prove basic statements.

There are two steps in the proof of the Fundamental Theorem of Arithmetic: (1) we first prove every integer at least two can be written as a product of prime powers, and then (2) prove that (up to reordering the factors, of course) there is only one way to do so. Thus, we first show existence and then uniqueness.

Existence is straightforward from the definition of primality. Recall an integer $n \ge 2$ is prime if it is only divisible by 1 and itself; 1 is called a unit, and all other positive integers are called composite. We declare 1 to be a unit and not a prime as otherwise unique factorization cannot hold; for example, $6$ $=$ $2 \cdot 3$ $=$ $1^{2011} \cdot 2 \cdot 3$. One way to prove the existence of a decomposition is by strong induction. Strong induction is similar to regular induction. Both first require us to show our statement $P(n)$ is true for $n=1$ (or another fixed integer $r$, depending on your counting). In normal induction one then shows that if $P(n)$ is true implies $P(n+1)$ is true, and concludes that $P(n)$ is true for all $n\geq r$. In strong induction, one shows that if $p(k)$ is true for all $k \le n$ implies $P(n+1)$ is true, and concludes that $P(n)$ is true for all $n\geq r$.\\

\begin{proof}[Proof of existence of factorization into primes] Clearly 2 is a product of prime powers, as 2 is prime. We now assume all integers $k$ with $2\le k \le n$ can be written as a product of prime powers. Consider the integer $n+1$. Either it is a prime (and we are done) or it is divisible by a prime. In the latter case, say $p$ divides $n+1$, which we write as $p|(n+1)$. We can then write $n+1$ as $pm$ for some $m < n$; by the induction assumption, $m$ can be written as a product of primes since it is at most $n$, completing the proof.
\end{proof}

We now turn to the more interesting part of the proof, namely the uniqueness. A common approach is to use the fact that if a prime $p$ divides a product $ab$, then either $p|a$ or $p|b$. Once we have this, the rest of the proof follows quickly.\\

\begin{proof}[Proof of uniqueness of factorization, given a prime divides a product implies the prime divides a factor] We proceed by contradiction. If unique factorization fails, then there must be some smallest integer, say $n$, where it first fails. Thus, we look at the smallest number with two distinct factorizations: $$n \ = \ p_1^{r_1} \cdots p_\ell^{r_\ell} \ = \ q_1^{s_1} \cdots q_\ell^{s_\ell}.$$ We claim that $p_k$ must divide both sides. It clearly divides the first product, and by our assumption it must divide the second. Why? We constantly group terms. We first write the second product as $$\left(q_1^{s_1} \cdots q_\ell^{s_\ell-1}\right) \cdot q_\ell.$$ Either $p_k|q_\ell$, or since we are assuming a prime cannot divide a product without dividing at least one factor, then $p_k$ must divide the term in parentheses. Arguing along these lines eventually leads us to $p_k$ dividing some $q_j$, and thus $p_k = q_j$. If we now look at $n/p_k$, we still have two distinct factorization; however, this contradicts the minimality of $n$, which we assumed was the smallest integer that had two factorizations. This completes the proof.
\end{proof}

We have thus reduced the proof of unique factorization to a prime divides a product if and only if it divides at least one term. As elementary as this sounds, this must be proved, especially as generalizations fail! For example, consider the ring $$\Z[\sqrt{-5}] \ = \ \{a + b \sqrt{-5}: a, b \in \Z\}.$$ In $\Z[\sqrt{-5}]$ we have $2\cdot 3 = (1+\sqrt{-5})(1-\sqrt{-5})$ and neither 2 nor 3 divide either factor on the right. For a more interesting example of a situation where unique factorization fails (though this one requires more background), consider the ring of trigonometric polynomials, namely functions of the form $$a_0 + \sum_{k=1}^n a_k \cos(kx) + \sum_{k=1}^n b_k \sin(kx)$$ where $a_i, b_i$ are real numbers and $n$ may be any non-negative integer (if $n=0$ the sum is considered vacuous). Hale Trotter \cite{Tr} observed that unique factorization fails here, as the Pythagorean formula implies $$\sin x \sin x \ = \ (1 - \cos x) \cdot (1 + \cos x),$$ and $\sin x$ divides neither $1-\cos x$ nor $1+\cos x$.

How can we prove a prime dividing a product divides at least one of the factors? The standard approach in many books is to use the Euclidean algorithm, or B\'ezout's identity, which follows from it; our conversations began in an attempt to see what was the minimum amount of machinery needed to prove this innocuous sounding claim. B\'ezout's identity (see for instance Chapter 1 of \cite{MT-B}) states that given any two integers $a$ and $b$ there are integers $x$ and $y$ such that $xa+yb = {\rm gcd}(a,b)$, where ${\rm gcd}(a,b)$ is the greatest common divisor of $a$ and $b$. Recall the greatest common divisor is defined as the largest integer which divides both $a$ and $b$.

How does this help us? Assume $p|ab$ but $p$ divides neither $a$ nor $b$. We then apply B\'ezout's identity to $p$ and $a$, noting that ${\rm gcd}(p,a) = 1$. Thus there are integers $x,y$ with $xp + ya = 1$. Multiplying both sides by $b$ we find $xpb + yab = b$. Since $p$ divides $p$ and $ab$ we see $p$ divides the left hand side, which implies $p$ divides the right hand side, namely $b$.

%\textbf{Add stuff about Euclidean algorithm.} It is possible to avoid the Euclidean algorithm by appealing to B$\acute{{\rm e}}$zout's identity, which states that if $a$ and $p$ are relatively prime then there are integers $x$ and $y$ such that $xa+yp = 1$. Multiplying both sides by $b$ gives $xab + ypb = b$, and as the right hand side is divisible by $p$ so too must $b$. More generally, one can always show that for any $a$ and $b$ there are integers $x$ and $y$ with $xa+yb = {\rm gcd}(a,b)$.

In discussing B\'ezout's identity, we were struck with how absurd it felt to have to spend time proving the following fact: if $d = {\rm gcd}(a,b)$ and $k|a$ and $k|b$ then $k|d$. \emph{Clearly} any factor of both $a$ and $b$ must divide their greatest common divisor, right? Unfortunately, this is not immediate from the definition of the gcd. The gcd is the largest number dividing both -- one must show every other divisor also divides this! Because of our intuition on the integers, it seems ridiculous to think that we could have $k|a, k|b$ without $k|{\rm gcd}(a,b)$, but we must prove this if we want to use the identity.

Fortunately, it turns out to be relatively simple to prove that every divisor divides the greatest common divisor, \emph{assuming we have a proof of unique factorization!}\\

\begin{proof}[Proof that $k|a$, $k|b$ implies $k|{\rm gcd}(a,b)$ (assuming unique factorization)] By unique factorization, there are unique decompositions of $a$ and $b$ as products of prime powers, say $$a \ = \ 2^{r_1} \cdot 3^{r_2} \cdots p_\ell^{r_\ell}, \ \ \ b \ = \ 2^{s_1} \cdot 3^{s_2} \cdots p_\ell^{s_\ell},$$ where $r_i, s_i$ are non-negative integers. Let $$t_i \ = \ \min(r_i, s_i), \ \ \ d \ = \ 2^{t_1} \cdot 3^{t_2} \cdots p_\ell^{t_\ell}.$$ We leave it as an exercise to the reader to show that $d$ is the greatest common divisor of $a$ and $b$. If now $k|a$ and $k|b$, then writing $$k \ = \ 2^{u_1} \cdot 3^{u_2} \cdots p_\ell^{u_\ell},$$ we see that $k|a$ implies $u_i \le r_i$ (the largest power of $p_i$ dividing $a$), and similarly $k|b$ implies $u_i \le s_i$. Thus $u_i \le \min(r_i,s_i) = t_i$, so $k|d$.
\end{proof}

We have just seen another application of unique factorization -- it allows us to justify using the name greatest common divisor! We thus return to the question of proving that if $p|ab$ then $p|a$ or $p|b$, as for the purpose of proving unique factorization we only need this, and not the far stronger statement that is B\'ezout's identity.

As the purpose of this article is to highlight techniques, we give two different proofs emphasizing different perspectives. In both proofs we do use one step which is similar to a key part of the proof of the Euclidean algorithm: if $a > p$ and $p$ does not divide $a$ then there is an integer $n \ge 1$ and an integer $r \in \{1,\dots,p-1\}$ such that $a=np+r$. To see this, just keep removing multiples of $p$ from $a$ until we are left with something between 1 and $p-1$, and since we are removing $p$ each time we cannot fall from a number exceeding $p$ to a number $0$ or less, as that would requiring removing at least $p+1$. It is important to note that in our arguments below, we are not using any results about the greatest common divisor, which is used in the standard proof of the Euclidean algorithm. Both proofs appeal to the following: if there is a set of non-negative integers having a certain property, then there is a smallest element of that set. \\

\begin{proof}[First proof that $p|ab$ implies $p|a$ or $p|b$.]

From our arguments before, we know it suffices to show that if $a$ and $p$ are relatively prime, then there exist integers $x$ and $y$ such that $xp + ya = 1$. (Remember we then multiplied by $b$ to find $xpb +yab = b$; if we assume $p|ab$ and $p$ does not divide $a$, then since $p$ divides the left hand side we find $p$ divides the right hand side, namely $b$.)

We prove a slightly more general version than we need: if $a$ and $b$ are relatively prime (which means the largest integer dividing both is 1), then there are integers $x$ and $y$ with $xa+yb=1$. We give a non-constructive proof that demonstrates the existence of $x$ and $y$.

First set
\[
A\ = \ \{xa+yb: x,y\in\Z\}\cap\ \{1,\ 2,\ 3,\ 4, \dots\}.
\] Thus $A$ is the set of positive integer values assumed as we vary $x$ and $y$ over all integers. As $a\in A$, clearly $A$ is non-empty and thus it has a smallest element, which we denote by $k$. This means $k$ can be written as $k=x_0a+y_0b$ for some integers $x_0,y_0$.

To complete the proof, we must show that $k=1$. We do this by first showing that $k$ divides
every element of $A$. Suppose that this were not the case. Then there would exist  a smallest element of $A$ that $k$ does not divide. As this element is in $A$, there must be integers $x_1$ and $y_1$ such that it equals $x_1a+y_1b$. Let $n$ be the largest natural number such that $nk < x_1a+y_1b<(n+1)k$. Why must there be such an $n$? The reason is we just keep adding multiples of $k$ until we first exceed our number $x_1a+y_1b$. It must be the case that
\[ 0\ < \ (x_1a+y_1b)-nk\ < \ k.  \] Why? If the middle expression equals 0 then we would have $x_1a+y_1b = qk$, which means $k$ divides $x_1a+y_1b$ and contradicts our assumption that $x_1a+y_1b$ is the smallest element in $A$ not divisible by $k$. A similar argument shows that it cannot equal $k$.

Since $k=x_0a + y_0b$, we see that $$x_1a + y_1b - k \ = \ (x_1-x_0)a + (y_1 - y_0)b.$$ As this number is positive and of the form $xa+yb$, we see it is in $A$. Further, it is strictly smaller than $k$, which contradicts $k$ being the smallest element of $A$. Thus our assumption that there is an element in $A$ not divisible by $k$ is false, and all elements of $A$ are divisible by $k$.

The proof is completed by looking at good choices for elements of $A$. If we take $x=1$ and $y=0$ we find $a\in A$, while taking $x=0$ and $y=1$ gives $b\in A$. Thus $k|a$ and $k|b$. As $a$ and $b$ have no common factors, there is only one possibility: $k=1$. Equivalently, our choice of $x_0$ and $y_0$ gives $$1 \ = \ x_0 a + y_0 b,$$ which is exactly what we wished to show.
\end{proof}

\begin{proof}[Second proof that $p|ab$ implies $p|a$ or $p|b$]

We proceed by contradiction.

\begin{enumerate}

\item Let $p$ be the smallest prime such that there is a product $ab$ with $p|ab$ but $p$ divides neither $a$ nor $b$.

\item Assuming such a prime $p$ exists, let $a$ and $b$ be such that the product $ab$ is smallest among all products of integers where $p$ divides the product but not the factors; if there are several decompositions giving the same smallest product, for definiteness take the one where $a$ is smallest.

\item We now show $a, b < p$. Clearly $a, b \neq p$. Assume $a > p$ (a similar argument holds if $b > p$). By our discussion above we can write $a = np + r$ with $0 < r < p$. This gives $ab = (np+r)b = nbp + rb$. As $p|ab$, subtracting gives $p|rb$; however, $rb < ab$, contradicting the minimality of the product $ab$. Thus $a, b < p$.

\item  As all integers have a prime decomposition, we write $$a \cdot b \ = \ q_{a,1} \cdots q_{a,\ell} \cdot q_{b,1} \cdots q_{b,m}.$$ As we are assuming $p|ab$, we may write $ab = np$, so $$np \ = \ q_{a,1} \cdots q_{a,\ell} \cdot q_{b,1} \cdots q_{b,m}.$$ We have assumed that $p$ is the smallest prime such that there is a product where $p$ divides the product but divides neither factor. As the factors of $a$ are at most $a<p$ (and those of $b$ are at most $b<p$), by induction each factor $q_{a,i}$ and $q_{b,j}$ divides either $n$ or $p$. As $p$ is prime, its only divisors are itself and 1. None of the divisors $q_{a,1}, \dots, q_{b,m}$ can divide $p$ as we have assumed $p$ divides neither $a$ nor $b$ (and if say $q_{a,1}|p$ then we would find $q_{a,1} = p$ and thus $p|a$).

\item Thus each factor divides $n$, so $$n \ = \ \widetilde{n} q_{a,1} \cdots q_{a,\ell} \cdot q_{b,1} \cdots q_{b,m} \ = \ \widetilde{n} ab.$$ As $ab = np$, we find $$n \ = \ n \widetilde{n} p,$$ which implies $\widetilde{n}p = 1$, which is impossible as $p \ge 2$. We have thus found a contradiction, and therefore there cannot be a smallest prime which divides a product without dividing at least one factor.
\end{enumerate}
\end{proof}

This completes our proof that a prime dividing a product must divide a factor, which is what we needed to prove unique factorization. Our hope in this note was twofold. The first was to highlight innocuous looking statements, stressing the need for careful proofs. It's very easy to be misled by notation. For example, the greatest common divisor of two numbers is actually just the greatest (i.e., largest) divisor of the numbers; we've inserted the word common, but we must justify its inclusion. The second was to see exactly how deep certain results are, specifically, exactly how much of certain results we need for our arguments.

%%%%%%%%%%%%%%%%%%%%%%%%%%%%%%%%%%%%%%%%%%%%%%%%%%%%%%%%%%%%%%%%%%%%%%%%%%%%%%%%%%%%%%%%%%%%%%%%%%%%%%%%%%%%

%%%%%%%%%%%%%%%%%%%%%%%%%%%%%%%%%%%%%%%%%%%%%%%%%%%%%%%%%%%%%%%%%%%%%%%%%%%%%%%%%%%%%%%%%%%%%%%%%%%%%%%%%%%%

\ \\


\begin{thebibliography}{BasBo2} % '2nd argument contains the widest acronym'

\bibitem[Cl]{Cl}
Clay Mathematics Institute, \texttt{http://www.claymath.org/millennium/Riemann\underline{\ }Hypothesis/}.

\bibitem[HW]{HW}
G. H. Hardy and E. Wright, \emph{An Introduction to the Theory of
Numbers}, 5th edition, Oxford Science Publications, Clarendon Press,
Oxford, $1995$.

\bibitem[MT-B]{MT-B}
S. J. Miller and R. Takloo-Bighash, \emph{An Invitation to Modern Number Theory}, Princeton University Press, Princeton, NJ, 2006, 503 pages.


\bibitem[Si]{Si}
C. E. Silva, \emph{Invitation to Real Analysis}. In preparation.

\bibitem[Tr]{Tr}
H. Trotter, \emph{An overlooked example of nonunique factorization}, Amer. Math. Monthly \textbf{95} (1988), no. 4, 339--342.

\bibitem[Wi]{Wi}
Wikipedia: Pages on B$\acute{{\rm e}}$zout's identity, Euclid's lemma, and the Fundamental Theorem of Arithmetic. Accessed October, 2010.



\end{thebibliography}
\end{document}